\documentclass[leqno,12pt]{amsart} 

\usepackage{amssymb}

\newtheorem{thm}{Theorem}[section]
\newtheorem{prop}[thm]{Proposition}

\newtheorem{cor}[thm]{Corollary}

\theoremstyle{definition}
\newtheorem{defn}[thm]{Definition}

\theoremstyle{remark}

\numberwithin{equation}{section}

\def\C{\mathbb{C}}

\def\R{\mathbb{R}}

\def\CP{\mathbb{CP}}

\def\E{\mathcal{E}}
\def\P{\mathbb{P}}

\def\D{\mathcal{D}}
\def\D{\mathcal{D}}

\newcommand{\de}{\partial}
\newcommand{\db}{\overline{\partial}}
\newcommand{\ddb}{{\partial }\overline{\partial}}

\begin{document}






\title[K\"ahler or balanced basis]{Holomorphic submersions onto K\"ahler or balanced manifolds}

\author{Lucia Alessandrini}
\address{ Dipartimento di Matematica\newline
Universit\`a degli Studi di Parma\newline
Parco Area delle Scienze 53/A\newline
I-43124 Parma
 Italy} \email{lucia.alessandrini@unipr.it}

\subjclass[2010]{Primary 53C55; Secondary 32J27, 32L05}


\keywords{K\"ahler manifolds,
balanced manifolds, sG manifolds, SKT manifolds, $p-$K\"ahler manifolds.}

\begin{abstract}
 We study many properties concerning weak K\"ahlerianity on  compact complex manifolds which admits a holomorphic submersion onto a K\"ahler or a balanced manifold. We get generalizations of some results of Harvey and Lawson (the K\"ahler case), Michelsohn (the balanced case), Popovici (the sG case) and others.

\end{abstract}

\maketitle

\section{Introduction}

It is well known that a compact holomorphic fibre bundle with K\"ahler basis and K\"ahler standard fibre does not carry, in general, a K\"ahler metric: this fact heavily depends on the cohomology of the total space, in particular on the vanishing of the cohomology class of the standard fibre. Simple examples are the Iwasawa manifold $I_3$, the Hopf manifolds and the Calabi-Eckmann spheres.

$I_3$ is a compact holomorphic fibre bundle on a two-dimensional complex torus $T_2$, whose standard fibre is a one-dimensional torus $T_1$ (see \cite{GH}, p. 444). $I_3$ is not K\"ahler because the homology class of the standard fibre vanishes (that is, the fibre bounds); nevertheless, $I_3$ is a balanced manifold.

Let us recall the definition of the Calabi-Eckmann spheres: $M_{u,v} := S^{2u+1} \times S^{2v+1}$, endowed with one of the complex structures of Calabi-Eckmann, is the total space of a (principal) holomorphic fibre bundle over the basis $\CP_u \times \CP_v$, with standard fibre (and structure group) a torus $T_1$ (in case $u=0$ or $v=0$, they are Hopf manifolds); $M_{u,v}$ is not K\"ahler nor balanced (see \cite{Mi}).

\medskip
We consider in the present paper two kinds of questions, namely:

i) We search suitable conditions which can be added to those on the basis, to get a K\"ahler or a balanced total space.

ii) If the basis is \lq\lq K\"ahler\rq\rq in a more general sense (i.e., it has a hermitian metric which is pluriclosed (SKT), or strongly Gauduchon, or hermitian symplectic \dots see section 2), we would like to get the same condition on the total space.

\medskip
As a matter of fact, we shall look at this kind of problems in a little more general setting, that is:

Let $M$ and $N$ be connected compact complex manifolds, with ${\rm dim} N = n > m = {\rm dim} M \geq 1$, and let $f : N \to M$ be a holomorphic submersion, where $a := n-m = {\rm dim} f^{-1} (x), \ x \in M,$ is the dimension of the standard fibre $F$.

Our hypotheses are of this kind:

a) $M$ has a K\"ahler or a balanced metric;

b) the class of the fibre $F$ does not vanish in a suitable cohomology group of $N$.
\medskip

We look for some \lq\lq $q-$K\"ahler\rq\rq properties on $N$: but before illustrating the results (collected in theorems 3.4, 3.5, 3.6, 3.8), we should explain precisely what are the right cohomology groups and what we mean with \lq\lq $q-$K\"ahler\rq\rq. 
This is not a simple matter at all, because almost everyone has given new names to the objects: we shall try to give also a  \lq\lq dictionary\rq\rq to understand the connection with other papers.
\medskip

Two old theorems can explain the background of our results, namely:

\begin{thm} {\rm (\cite{HL}, Theorem 17)}
Suppose $f:N \to M$ is a holomorphic submersion with $1-$dimensional fibres onto a K\"ahler manifold $M$. Then there exists a K\"ahler metric on $N$ if and only if the fibre of $f$ is not a $(1,1)-$component of a boundary.
\end{thm}

\begin{thm} {\rm (\cite{Mi}, Theorem 5.5)}
Suppose $f:N \to C$ is a holomorphic map from a compact complex manifold onto a curve $C$. Then there exists a balanced metric on $N$ if no positive linear combination of irreducible components of fibres
 of $f$ is an $(n-1,n-1)-$component of a boundary, and the non-singular fibres of $f$ are balanced.
\end{thm}
\medskip

We refer to our paper  \cite{A11} for the full generality: here we recall only the basic definitions, starting from the cases $p=1$ and $p=n-1$, which are principally involved in our present results.

\section{Preliminaries}\label{S2}

Let $N$ be a compact complex manifold of dimension $n \geq 2$, let $p$ be an integer, $1 \leq p \leq n-1$. 
As  regards forms and currents, we shall use mainly the notation of \cite{De}. 

A $(k,k)$-current $T$  is a current of
bidegree $(k,k)$ or bidimension $(p,p)$, where $p+k=n$; 
$T \in {\mathcal{D}}_{p,p}'(N)_{\R}$ means that $T$ is a
real $(k,k)$-current on $N$; in particular, if $T$ is a positive $(k,k)$-current 
($T \geq 0$), then it is real.

\medskip

We shall need de Rham cohomology, and also Bott-Chern and Aeppli cohomology (for which the notation is not standard): both of them can
be described using forms or currents of the same bidegree: 

$$H_\R ^{k,k}(N) :=\frac{\{ \varphi \in {\E}^{k,k}(N)_\R;
d\varphi =0\}}{\{d\psi ;\psi \in {\E}^{2k-1}(N)_\R\}}\simeq\frac{\{T \in
{\D}_{p,p}'(N)_{\R}; dT =0\}}{\{dS ; S \in  {\D}_{2p+1}'(N)_{\R}
\}}.$$
$$H_{\ddb}^{k,k}(N) = \Lambda_\R ^{k,k}(N)= H_{BC}^{k,k}(N) :=\frac{\{ \varphi \in {\E}^{k,k}(N)_\R;
d\varphi =0\}}{\{i\partial\overline{\partial}\psi ;\psi \in {\E}^{k-1,k-1}(N)_\R\}}$$ $$\simeq\frac{\{T \in
{\D}_{p,p}'(N)_{\R}; dT =0\}}{\{i\partial\overline{\partial}A ; A \in  {\D}_{p+1,p+1}'(N)_{\R}
\}}.$$
$$H_{\de + \db}^{k,k}(N) =V_\R ^{k,k}(N)= H_{A}^{k,k}(N) :=\frac{\{ \varphi \in {\E}^{k,k}(N)_\R;
i\ddb\varphi =0\}}{\{\varphi = \de \overline\eta + \db \eta ; \eta \in {\E}^{k,k-1}(N)\}}$$ $$\simeq\frac{\{T
\in {\D}_{p,p}'(N)_{\R}; i\ddb T =0\}}{\{ \de \overline S + \db S ; S \in  {\D}_{p,p+1}'(N)
\}}.$$

\bigskip
In general when the class of a current vanishes in one of the previous cohomology groups, we say that the current {\bf \lq\lq bounds\rq\rq}.

We collect definitions and characterization's results in the following definition (see \cite{A11}).

\medskip

\begin{defn} 

\begin{enumerate}

\item {\bf Characterization of  $p-$K\"ahler (pK) manifolds.}

$N$ has a strictly weakly positive (i.e. transverse) $(p,p)-$form $\Omega$ with $\de \Omega = 0$,  if and only if $N$ has no strongly positive currents $T \neq 0$, of bidimension $(p,p)$, such that $T = \de  \overline S + \db S$ for some current $S$ of bidimension $(p,p+1)$ (i.e.  $T$  \lq\lq bounds\rq\rq in $H_{\de + \db}^{k,k}(N)$, i.e. $T$ is the $(p,p)-$component of a boundary).

\item {\bf Characterization of  weakly $p-$K\"ahler (pWK) manifolds.}

$N$ has a strictly weakly positive (i.e. transverse) $(p,p)-$form $\Omega$ with $\de \Omega = \ddb \alpha$ for some form $\alpha$,  if and only if $N$ has no strongly  positive currents $T \neq 0$, of bidimension $(p,p)$, such that $T = \de  \overline S + \db S$ for some current $S$ of bidimension $(p,p+1)$ with $\ddb S = 0$ (i.e.  $T$  is closed and \lq\lq bounds\rq\rq in $H_{\de + \db}^{k,k}(N)$). 

\item {\bf Characterization of  $p-$symplectic (pS) manifolds.}

$N$ has a real $2p-$form $\Psi = \sum_{a+b=2p} \Psi^{a,b}$, such that $d \Psi = 0$ and the 
$(p,p)-$form $\Omega := \Psi^{p,p}$ is 
strictly weakly positive,   if and only if $N$ has no strongly  positive currents $T \neq 0$, of bidimension $(p,p)$, such that $T = d S$ for some current $S$  (i.e.  $T$  is a boundary in de Rham cohomology). 

\item {\bf Characterization of   $p-$pluriclosed (pPL) manifolds.}

$N$ has a strictly weakly positive $(p,p)-$form $\Omega$ with $\ddb \Omega = 0$,  if and only if $N$ has no strongly  positive currents $T \neq 0$, of bidimension $(p,p)$, such that $T =i   \ddb A$ for some current $A$ of bidimension $(p+1,p+1)$ (i.e.  $T$  \lq\lq bounds\rq\rq in $H_{\ddb}^{k,k}(N)$). 
\end{enumerate}

 \end{defn}
\medskip

{\bf 2.2 Remark.} The technique used to prove the previous characterization statements stems from the work of Sullivan \cite{Su}, and is based on the Hahn-Banach Separation Theorem (on dual spaces of forms and currents): see \cite{A11} for the proofs. 

\medskip
{\bf 2.3 Remark.} In particular, notice that the currents which are involved are positive in the sense of Lelong, i.e. strongly positive, so that the dual cone is that of weakly positive forms. To be precise, we should define  weakly positive, positive, strongly positive currents (see  \cite{HK}, \cite{A11}), but the wider class, that of weakly positive currents, is enough for our purpose, hence we speak of  {\it positive} currents in general.

\medskip
{\bf 2.4 Remark. } As regards Definition 2.1(3), let us write the condition $d \Psi = 0$ in terms of a condition on $\de \Omega$, as in the other statements; $d \Psi = 0$ is equivalent to:

i) $\db \Psi^{n-j,2p-n+j} + \de \Psi^{n-j-1,2p-n+j+1}=0$, for $j=0, \dots, n-p-1$, when $n\leq 2p$

and

ii) $\de \Psi^{2p,0}=0, \ \db \Psi^{2p-j,j} + \de \Psi^{2p-j-1,j+1}=0$, for $j=0, \dots, p-1$, when $n > 2p.$

In particular, $\de \Omega = \de \Psi^{p,p}= - \db \Psi^{p+1,p-1}$ (which is the sole condition when $p=n-1$).

\bigskip
When $M$ satisfies one of the characterization theorems given in Definition 2.1, in the rest of the paper we will call it generically a {\bf \lq\lq$p-$K\"ahler\rq\rq manifold}; the form $\Omega$ is said to be {\bf \lq\lq closed\rq\rq}.
Notice also that:
$pK \Longrightarrow pWK  \Longrightarrow pS  \Longrightarrow pPL.$ 

As regards examples and differences among these classes of manifolds, see \cite{A11}: $p-$K\"ahler and $p-$symplectic manifolds had been defined in \cite{AA}.

\bigskip

{\bf 2.5 The case $p=1$.} For $p=1$, a transverse form is the fundamental form of a hermitian metric, so that we can speak of $1-$K\"ahler, weakly $1-$K\"ahler,  $1-$symplectic, $1-$pluriclosed {\it metrics}. 

Notice that, while  a $1-$K\"ahler manifold is simply a K\"ahler manifold, the $1-$symplec\-tic condition means that there is a symplectic $2-$form $\Psi$ which tames the given complex structure $J$ (in the sense of Mc Duff and Gromov, i.e. $\Psi_x(v, Jv) > 0, \ \forall \ v \in T_xM$, 
see
\cite{MS}, \cite{Gr}; see moreover \cite{Su}, pp. 249-252); we get a hermitian metric with fundamental form $\alpha$ (not closed, in general). $1-$symplectic manifolds are also called {\it holomorphically tamed}, or
{\it hermitian symplectic} (\cite{ST}).
In \cite{E}, pluriclosed (i.e. $1-$pluriclosed) metrics are defined (see also \cite{ST}), while in  \cite{FT} a 1PL metric (manifold) is called a {\it strong K\"ahler metric (manifold) with torsion} (SKT). 

\bigskip
{\bf 2.6 The case $p=n-1$.} For  $p = n-1$, we get a hermitian metric too, because  every transverse $(n-1,n-1)-$form $\Omega$ is in fact given by $\Omega = \omega^{n-1}$, where $\omega$ is a transverse $(1,1)-$form (see f.i. \cite{Mi}, p. 279). 

This case was studied by Michelsohn  in \cite{Mi}, where $(n-1)-$K\"ahler manifolds are called {\it balanced} manifolds. 

Moreover, $(n-1)-$symplectic manifolds are called {\it strongly Gauduchon manifolds (sG)} by Popovici (compare Remark 2.4 and Definition 2.1(3) with  \cite{Po1}, Definition 4.1 and Propositions 4.2 and 4.3; see also \cite{Po2}), while $(n-1)-$pluriclosed metrics are called {\it standard} or {\it Gauduchon metrics}. Recently, weakly $(n-1)-$K\"ahler manifolds have been called {\it superstrong Gauduchon (super sG)} (\cite{PU}).

\medskip
{\bf 2.7 Remark.} Every compact complex manifold supports Gauduchon metrics: in fact, by the characterization in Definition 2.1(4), if 
$T$ is a positive $(1,1)-$current, such that $T = \ddb A$, $A$ turns out to be a plurisubharmonic function; but $N$ is compact, so that $A$ is constant, and $T=0$.
\medskip

{\bf  2.8 Remark.} As regards compact complex surfaces ($n=2$), we have: 

Every surface is 1PL (SKT), because $1=n-1$; moreover, there is only a class of special surfaces, those which are K\"ahler (i.e. balanced), because (see \cite{La}):
$$1K \iff b_1 {\rm \  is \ even}  
\iff 1S.$$

The Hopf surface is not in this class.

Let us notice that this regards manifolds, but not metrics, as it involves the non-existence of currents! \bigskip

{\bf 2.9 The case $1 < p < n-1$.} When $1 < p < n-1$, and $\omega$ is a transverse $(1,1)-$form, $d \omega^p = 0$ implies $d \omega = 0$; moreover, a transverse $(p,p)-$form $\Omega$ is not necessarily of the form  $\Omega = \omega^{p}$, where $\omega$ is a transverse $(1,1)-$form (see also section 4). 

Hence in the intermediate cases ($1 < p < n-1$) the $(p,p)-$form $\Omega$ in Definition 2.1 is not of the form $\Omega = \omega^{p}$, in general. Therefore we will not look for \lq\lq good\rq\rq hermitian metrics, but will instead handle transverse forms or positive currents, as done in Definition 2.1.
\medskip

After all, let us recall a very useful result:

{\bf The division theorem} (see \cite {Le}, Theorem 2, p. 69).

{\it Let $\psi$ be a positive $(1,1)-$form of rank $m$ on a manifold $N$ (i.e. $\psi^m \neq 0, \psi^{m+1} =0$), and let $t$ be a positive current on $N$ of bidegree $(q,q)$, such that $t \wedge \psi =0$.

\begin{enumerate}
\item If $m>q$, then $t=0$.

\item If $m \leq q$, then there is a unique positive current $R$ of bidegree $(q-m,q-m)$ on $N$ such that $t=R \wedge \psi^m$. In particular, if $q=m$, there is a positive measure $\mu$ on $N$ such that $t= \mu \psi^m$.
\end{enumerate}}

\medskip

\section{Results}\label{S3}

Let $M$ and $N$ be connected compact complex manifolds, with ${\rm dim} N = n > m = {\rm dim} M \geq 1$, and let $f : N \to M$ be a holomorphic submersion, where $a := n-m = {\rm dim} f^{-1} (x), \ x \in M,$ is the dimension of the standard fibre $F$.
\medskip

As regards the push forward of a \lq\lq$p-$K\"ahler\rq\rq property, we have:

\begin{prop}
Let $f : N \to M$ be as above. If $N$ is \lq\lq$p-$K\"ahler\rq\rq for some $p, \ a < p \leq n-1$, then $M$ is \lq\lq$(p-a)-$K\"ahler\rq\rq. In particular, if $N$ is balanced, then $M$ is balanced too.
\end{prop}

{\bf Proof.} If $\Omega$ is a  \lq\lq closed\rq\rq transverse $(p,p)-$form on $N$, then $f_*\Omega$ is a \lq\lq closed\rq\rq transverse $(p-a,p-a)-$form on $M$.

\medskip
A deeper result is due to Varouchas (see \cite{V}):
\begin{thm}
Let $f : N \to M$ be a surjective holomorphic map with equidimensional fibres. If $N$ is K\"ahler, then $M$ is K\"ahler too.
\end{thm}
\medskip

Suppose on the contrary that $M$ has a K\"ahler or a balanced metric, with fundamental form $\omega$; our aim is to prove that $N$ is \lq\lq$p-$K\"ahler\rq\rq for some $p$; but pulling back $\omega$ we get the $(1,1)-$form $f^*\omega$ on $N$, which is no more strictly positive, but only $f^*\omega \geq 0$. Thus we switch to currents, and try to prove that there are no positive currents on $N$ which \lq\lq bound\rq\rq, as said in the characterization theorems (see Definition 2.1). For brevity, we shall study all cases together: this choice may make the following statements dull reading, but we discuss each case separately after the proofs.
\medskip

Fix an index $p, \ 1 \leq p \leq n-1$: in order to apply the division theorem, choose a \lq\lq bad\rq\rq  current $T$ on $N$, i.e. a positive current $T$ of bidimension $(p,p)$ with $T = \de  \overline S + \db S$ for some current $S$ of bidimension $(p,p+1)$ as in Definition 2.1, or $T = i\ddb A$ for some current $A$ of bidimension $(p+1,p+1)$; the aim is to conclude that $T=0$. 
\medskip

Consider 
$T \wedge f^*\omega ^h, \ 1 \leq h \leq min \{m,p\}$.

\bigskip
{\bf Step 1.} In the previous notation, if $d \omega^h=0$, then $T \wedge f^*\omega ^h$ is also \lq\lq bad\rq\rq.

{\bf Proof of Step 1.} Suppose $\de \omega^h=0$. Then if $T = \de  \overline S + \db S$, we get 
$$ \de ( \overline S \wedge f^*\omega ^h) + \db (S \wedge f^*\omega ^h) = \de \overline S \wedge f^*\omega ^h + \db S \wedge f^*\omega ^h = T \wedge f^*\omega ^h,$$

with $ \de ( S \wedge f^*\omega ^h) = \de S \wedge f^*\omega ^h$ and 
$ \ddb ( S \wedge f^*\omega ^h) = \ddb S \wedge f^*\omega ^h$;
thus we have on $T \wedge f^*\omega ^h$ the same conditions as on $T$.

If $T = i\ddb A$, we get $T \wedge f^*\omega ^h = i \ddb (A \wedge f^*\omega ^h)$.

\bigskip

To use the division theorem, we need $T \wedge f^*\omega ^h=0$:

\medskip
{\bf Step 2.} In the previous notation,  suppose $d \omega^h=0$. Then  $T \wedge f^*\omega ^h=0$ in the following cases:

\begin{enumerate}
\item $p=h$.

\item $p>h$ and $N$ is \lq\lq$(p-h)-$K\"ahler\rq\rq.

\end{enumerate}

{\bf Proof of Step 2.} \begin{enumerate}
\item When $p=h$, the current $T \wedge f^*\omega ^h$ has maximum degree, so that  $T \wedge f^*\omega ^h= \mu dV,$ where $dV$ is a volume form on $N$ and $\mu$ is a positive measure on $N$. But $\int_N \mu dV =0$, because $T \wedge f^*\omega ^h$  \lq\lq bounds\rq\rq (Step 1) and $N$ is compact, hence $\mu =0$.

\item When $p>h$, by Step 1, $T \wedge f^*\omega ^h$ is a \lq\lq bad\rq\rq current  of bidimension $(p-h,p-h)$ on a \lq\lq$(p-h)-$K\"ahler\rq\rq manifold, thus it vanishes.

\end{enumerate}

\bigskip
{\bf Step 3.} Let us apply now the division theorem with $\psi = f^*\omega$ ($rk \psi = m$), and with $t = T \wedge f^*\omega ^{h-1}$
($t=T$ in case $h=1$): this assures $t \wedge \psi =T \wedge f^*\omega ^h$. We get:

(i) If $T \wedge f^*\omega ^h=0$ and $a < p-h+1$, then $t=0$.

(ii) If $T \wedge f^*\omega ^h=0$ and $a = p-h+1$, then there exists a positive measure $\mu$ on $N$ such that $t= \mu f^*\omega ^{m}$.

{\bf Proof of Step 3.} 

(i) We get $m > q$, where $q$ is the bidegree of $t$, since $q=n-p+h-1$, but $a < p-h+1$; thus by the division theorem, $t=0$.

(ii) We have only to check, as before,  that $q=m$.

\bigskip

Recall that our goal is $T=0$.

\medskip
{\bf Step 4.} In case (i) ($a < p-h+1$ and $T \wedge f^*\omega ^h=0$), we get precisely $T=0$.

{\bf Proof of Step 4.} Obvious when $h=1$; in general, we get $T \wedge f^*\omega ^{h-1}=0$, thus we can apply the division theorem again, using $T \wedge f^*\omega ^{h-2}$ and getting $T \wedge f^*\omega ^{h-2}=0$, and so on, until  $T=0$.

\bigskip
{\bf Step 5.} In case (ii) ($a = p-h+1$ and $T \wedge f^*\omega ^h=0$), if moreover $\ddb t = \ddb (T \wedge f^*\omega ^{h-1}) = \ddb (\mu f^*\omega ^{m}) =0$,
then there exists a positive measure $\nu$ on $M$ such that $\mu=  f^*\nu$, so that $t =  f^*(\nu \omega ^{m})$.

{\bf Proof of Step 5.} The proof goes as in Lemma 18 in \cite{HL}: \lq\lq Suppose $f:X \to Y$ is a holomorphic submersion with one-dimensional fibres, and suppose $t$ is a positive current of bidimension $(1,1)$ on $X$. Then the push-forward $f_*t$ of $t$ to $Y$ is zero if and only if $t= ||t|| \overline F$, where $\overline F$ is the field of unit 2-vectors tangent to the fibre. If, in addition, $t$ satisfies the equation $\ddb t = 0$, then $t =f^*\nu$, for some non-negative density $\nu$ on $Y$\rq\rq.

Notice that the analogous of this Lemma when $a > 1$ is no more true, but it is not hard to prove that in our hypotheses the second part of the Lemma also holds when $a \neq 1$, because for dimensional reasons $\de (f^*\omega ^{m}) = f^*(\de \omega ^{m}) = 0$, thus 
$0 = \ddb (\mu f^*\omega ^{m}) = \ddb \mu \wedge f^*\omega ^{m}$.
This implies that, in the fibre directions, the measure $\mu$ is harmonic; since the fibres are compact, we conclude that $\mu$ is independent on fibre coordinates, i.e., there exists a positive measure $\nu$ on $M$ such that $\mu=  f^*\nu$.

\bigskip
 
 We get finally the following Proposition:
 
\begin{prop} In the above notation, suppose $T \wedge f^*\omega ^h=0$; we get $T=0$ when:

\begin{enumerate}
\item  $h=1$ and $p>a$;

\item $h=1$, $p=a$ and moreover the generic fibre $F$ does not \lq\lq bound\rq\rq in $N$;

\item $m>1$,  $h=m-1$, and $p=n-1$ (thus $p-h+1>a$).

\end{enumerate}

\end{prop}

{\bf Proof.} (1) and (3) are proved by Step 4.

As regards (2),  it holds $T= t =  \mu f^*\omega ^{m}$, because we are in case (ii) of Step 3. Notice that $\ddb T = 0$ since $T$ is \lq\lq bad\rq\rq, then by Step 5 there exists a positive measure $\nu$ on $M$ such that $\mu=  f^*\nu$, i.e. $T =  f^*(\nu \omega ^{m})$.
 
 For every $x \in M$, put $c := \int_M \nu \omega^m$. Then $\{\nu \omega^m\} = c \{ \delta_x \omega^m\}$ as homology classes in $M$, since the homology is one-dimensional in top degree.
 
 Pulling back by $f$, we have $c \{f^{-1}(x)\} = \{T\} =0$, but the generic fibre $F$ does not \lq\lq bound\rq\rq in $N$, hence $c=0$, so that $T=0$.

\bigskip

{\bf Claim.} Since the cohomology class of every fibre of a holomorphic submersion is the same,  in our setting we can consider the following {\it homological conditions} on $N$, (which does not depend on the index $p$):
$$(HC)_K \Longrightarrow (HC)_{WK} \Longrightarrow (HC)_S \Longrightarrow (HC)_{PL},$$  
where 

{\bf (\lq\lq HC\rq\rq)}: the generic fibre $F$ of $f: N \to M$ does not \lq\lq bound\rq\rq in $N$.
\medskip

It is clear that when $N$ is \lq\lq $a-$K\"ahler\rq\rq then (\lq\lq HC\rq\rq) holds; moreover, since the current given by the integration on $F$ is a {\it closed}  positive current of bidimension $(a,a)$ on $N$,  $(HC)_K=(HC)_{WK}$. 
\bigskip

Thus we got $T=0$ in all cases, so that $N$ is \lq\lq $p-$K\"ahler\rq\rq: let us collect our results in the following theorems, starting from low dimensional manifolds.
\medskip

\begin{thm} Let $M$ and $N$ be compact complex manifolds, with ${\rm dim} N = n > m = {\rm dim} M =1$, and let $f : N \to M$ be a holomorphic submersion, where $a := n-1 = {\rm dim} f^{-1} (x), \ x \in M,$ is the dimension of the standard fibre $F$.

\begin{enumerate}
\item If $n=2$, then: $N$ is \lq\lq$1-$K\"ahler\rq\rq if and only if (\lq\lq HC\rq\rq) holds.

\item If $n>2$, and $N$ is \lq\lq$(n-2)-$K\"ahler\rq\rq, then: $N$ is \lq\lq$(n-1)-$K\"ahler\rq\rq if and only if (\lq\lq HC\rq\rq) holds.

\end{enumerate}

\end{thm}

{\bf Proof.} It is a particular case of Theorem 3.6.

\medskip
{\bf Remarks on Theorem 3.4.} 
\begin{enumerate}
\item The case PL is not significative, since every compact manifold is $(n-1)$PL.

\item If $N$ is a surface, all \lq\lq K\"ahler\rq\rq conditions are equivalent, except PL (see Remark 2.8):  thus the results we got are nothing but Theorem 17  in \cite{HL} (see also \cite{Mi}, Corollary 5.8).

\item Theorem 3.4(2) in case K is in fact a particular case of Theorem 5.5 in \cite{Mi} (see Theorem 1.2), because when  $N$ is $(n-2)$K, then every fibre  is balanced (pulling back the form from $N$ to every fibre). Cases WK and S seem to be new.

\end{enumerate}

\bigskip
\begin{thm} Let $M$ and $N$ be compact complex manifolds, with ${\rm dim} N = n > m = {\rm dim} M =2$, and let $f : N \to M$ be a holomorphic submersion, where $a := n-2 = {\rm dim} f^{-1} (x), \ x \in M,$ is the dimension of the standard fibre $F$. Suppose $M$ is K\"ahler, (i.e. balanced, 1S, 1WK).

\begin{enumerate}
\item If $N$ is \lq\lq$(n-2)-$K\"ahler\rq\rq, then it is also \lq\lq$(n-1)-$K\"ahler\rq\rq.

\item If $n=3$, then $N$ is \lq\lq$1-$K\"ahler\rq\rq if and only if (\lq\lq HC\rq\rq) holds.

\item If $n>3$, and $N$ is \lq\lq$(n-3)-$K\"ahler\rq\rq, then: $N$ is \lq\lq$(n-2)-$K\"ahler\rq\rq if and only if (\lq\lq HC\rq\rq) holds.

\end{enumerate}

\end{thm}

{\bf Proof.} It is a particular case of Theorem 3.6.

\medskip
{\bf Remarks on Theorem 3.5} (See also Remarks on Theorem 3.6).
\begin{enumerate}
\item In Theorem 3.5(1), the case PL is not significative, since every compact manifold is $(n-1)$PL. 

\item Since $(HC)_K=(HC)_{WK}$, if $n=3$ we get that $N$ is 1WK if and only if it is K\"ahler. 

\item If $n=3$, compare Theorem 3.5(2), case K, with Theorem 17  in \cite{HL} (i.e. Theorem 1.1 here). 

\item Consider the fibration $I_3 \to T_2$ (see Section 1), and recall that on $I_3$, all \lq\lq$p-$K\"ahler\rq\rq 

\noindent conditions are equivalent, since it is holomorphically parallelizable (\cite{AB}). This example shows that $(HC)_K$ is not a necessary condition to be balanced: in fact $I_3$ is \lq\lq$(n-1)-$K\"ahler\rq\rq but (\lq\lq HC\rq\rq) does not hold.

\end{enumerate}

\medskip

\begin{thm} Let $M$ and $N$ be compact complex manifolds, with ${\rm dim} N = n > m = {\rm dim} M \geq 3$, and let $f : N \to M$ be a holomorphic submersion, where $a := n-m = {\rm dim} f^{-1} (x), \ x \in M,$ is the dimension of the standard fibre. 
Suppose $M$ is K\"ahler.

\begin{enumerate}
\item If $N$ is \lq\lq$p-$K\"ahler\rq\rq, with $p \geq a$, then it is also \lq\lq$(p+1)-$K\"ahler\rq\rq, \dots, \lq\lq$(n-1)-$K\"ahler\rq\rq.

\item If $a=1$, then $N$ is \lq\lq$1-$K\"ahler\rq\rq if and only if (\lq\lq HC\rq\rq) holds.

\item If $a>1$, and $N$ is \lq\lq$(a-1)-$K\"ahler\rq\rq, then: $N$ is \lq\lq$a-$K\"ahler\rq\rq if and only if (\lq\lq HC\rq\rq) holds.

\end{enumerate}

\end{thm}
{\bf Proof.} As regards (1), since we can take $h=1$ ($\de \omega = 0$),  we get  from Proposition 3.3(1) together Step 2.(2) that $N$ is \lq\lq$(p+1)-$K\"ahler\rq\rq (put $q =p+1$, so that  the \lq\lq bad\rq\rq current $T$ has bidimension $(q,q)$; obviously it holds $q > a \geq 1$). Then we may apply the same over and over.

Moreover, we get (2) from Proposition 3.3(2) together Step 2.(1), and we get (3) from Proposition 3.3(2) together Step 2.(2) (notice that in both statements, one side is straightforward).

\medskip
{\bf Remarks on Theorem 3.6.} 
\begin{enumerate}

\item Since $(HC)_K=(HC)_{WK}$, if $a=1$ we get that $N$ is 1WK if and only if it is K\"ahler. 

\item Theorem 3.6(2) for the case K is exactly Theorem 17  in \cite{HL} (Theorem 1.1 here). In that paper, Harvey and Lawson asked for the case when the non-K\"ahler property can be characterized by holomorphic chains: this is the case here.
Theorem 3.6(2) for the case PL was proved  in \cite{E} (Theorems 4.5 and 4.6).

\item Moreover,  when $a=1$ and (\lq\lq HC\rq\rq) holds, $N$ is also \lq\lq$2-$K\"ahler\rq\rq, \dots, \lq\lq$(n-1)-$K\"ahler\rq\rq thanks to Theorem 3.6(1) (this is not obvious, compare Section 4). 

\item Theorem 3.6(3) in case K is in fact a generalization of Theorem 1.2 (when $f$ is a holomorphic submersion), because when  $N$ is $(a-1)$K, then every fibre  is balanced; nevertheless, we get that $N$ is $aK$, not only balanced. The other cases (WK, S, PL) seem to be completely new.

\end{enumerate}

\medskip

\begin{cor} Let $f : N \to M$ be a holomorphic submersion as above, with ${\rm dim} N = n =m+1, \ m = {\rm dim} M \geq 2$.
Then $N$ is K\"ahler if and only if $M$ is K\"ahler and $(HC)_K$ holds.
\end{cor}

{\bf Proof.} Use Theorem 3.6(2) and Theorem 3.2 (this result is obvious in the present situation).

\medskip

\begin{thm} Let $M$ and $N$ be compact complex manifolds, with ${\rm dim} N = n > m = {\rm dim} M \geq 3$, and let $f : N \to M$ be a holomorphic submersion, where $a := n-m = {\rm dim} f^{-1} (x)$ is the dimension of the standard fibre $F$. Suppose $M$ is balanced. Then if $N$ is \lq\lq$a-$K\"ahler\rq\rq,  it is also \lq\lq$(n-1)-$K\"ahler\rq\rq.
\end{thm}

{\bf Proof.} Since here we can take $h=m-1$,  we get the result from Proposition 3.3(3) together with Step 2.(2).

\bigskip
{\bf 3.9 Remark.} Corollary 3.7 may suggest the following conjecture: 

\lq\lq Also when $a > 1$ (and $M$ is K\"ahler), the condition (\lq\lq HC\rq\rq) is sufficient to give to $N$ some \lq\lq$p-$K\"ahler\rq\rq structure.\rq\rq

Taking in account the statements 3.6(1) and 3.6(3), the most natural property to look for is an  \lq\lq$(a+1)-$K\"ahler\rq\rq structure on $N$: this fact would imply, by Theorem 3.6(1), that $N$ will be also \lq\lq$(a+2)-$K\"ahler\rq\rq, \dots, \lq\lq$(n-1)-$K\"ahler\rq\rq. 

Looking at Proposition 3.3(1), we have only to assure $T \wedge f^*\omega =0$, starting from $d \omega =0$ on $M$, and from a \lq\lq bad\rq\rq  current $T$ on $N$. This was done in \cite{Mi}, proof of Theorem 5.5, when $M$ is not only K\"ahler, but it is a curve. Nevertheless, that highly non-trivial construction does not work in our case, because the forms $\varphi_{\varepsilon}$, whose role is that to squeeze the support of the current $T \wedge f^*\omega$ in small tubular neighborhoods of the fibres, are no more closed and positive.

More than that, we give here an example to show that the above conjecture cannot hold. Let us consider the compact complex threefold $M_{1,1}$ given in Section 1; $M_{1,1}$ is not 2S = $(n-1)$S (i.e. strongly Gauduchon) since, by definition, it has non-trivial hypersurfaces (thus closed positive currents of bidimension (1,1)), but all them bound in the homology, since $H_2(M_{1,1}) =0$. Hence by Definition 2.1(3) this manifold is not 2S (and therefore not 1S, by Theorem 3.5(1)).

Let us consider $N = M_{1,1} \times \C\P_2$ and the holomorphic submersion given by the projection $\pi : N \to \C\P_2$.
By the K\"unneth formula, in the homology of $N = M_{1,1} \times \C\P_2$, the class of the fibre $M_{1,1}$ does not vanish, thus the condition $(HC)_S$ holds.

But $N = M_{1,1} \times \C\P_2$ has no $p-$symplectic structure at all:

not for $p=1, p=2$ since it contains $M_{1,1}$ as a submanifold,

not for $p=3, p=4$ due to Proposition 3.1 applied to the projection $f: N \to M_{1,1}$.

This proves that the sole condition \lq\lq (HC)\rq\rq does not assure a \lq\lq$p-$K\"ahler\rq\rq structure, when the dimension of the fibre is bigger than one.
The {\it correct} statement is indeed Theorem 4.1(ii).
\bigskip

\section{Conclusions and applications}\label{S4}

Let $N$ be a complex manifold of dimension $n \geq 3$, let $p$ be an integer, $1 \leq p \leq n-1$.
In section 2, for $p < n-1$, we defined $p-$K\"ahler manifolds not by means of a hermitian metric, but only using a strictly weakly positive $(p,p)-$form (notice also that, by our choice of a strictly {\it weakly} positive $(p,p)-$form $\Omega$, we cannot deduce that $\Omega \wedge \Omega$ is a $2p-$K\"ahler form): the basic motivation stems from the following observation (which can be directly checked):
 for $p < n-1$, when $d \omega^p =0$, then $d \omega =0$; thus this kind of $p-$K\"ahler manifolds, where $\Omega = \omega^p$, are simply the K\"ahler manifolds.
 \medskip
 
 On the other hand, since $\de \omega^p = p \omega^{p-1} \wedge \de \omega$, $\de \omega =0$ implies $\de \omega^p =0$, so that a K\"ahler manifold is $pK$ for all $p$:
 this does not work in the pWK case and in the pPL case, because
 
 $\de \omega^p = p \omega^{p-1} \wedge \de \omega = p \omega^{p-1} \wedge \ddb \alpha$, whereas we need  $\de \omega^p = \ddb \beta$; moreover,
 
 $\ddb \omega^p = p \omega^{p-2} \wedge ((p-1) \de \omega \wedge \overline{ \de \omega} + \omega \wedge \ddb \omega)$, and in particular 
 $\ddb \omega^2 = 2 ( \de \omega \wedge \overline{ \de \omega} + \omega \wedge \ddb \omega)$.
 
 Thus only  $\ddb \omega^2 =0, \ddb \omega=0$ implies $\ddb \omega^p =0$ for every $p$.
 \medskip
 
 In case PL, this kind of metrics / manifolds have been considered: in particular, the following conditions were studied, on a strictly positive $(1,1)-$form $\omega$ on an $n-$dimensional manifold $N$, $n \geq 3$.
 
 \begin{enumerate}
\item $\ddb \omega^k = 0$ ($k-$SKT, see \cite{IP};  $k=1$ corresponds to SKT, i.e. 1PL, $k=n-2$ is the {\it astheno-K\"ahler} condition 
 of Yost and Yau, $k=n-1$ corresponds to Gauduchon metrics)
 
\item $\ddb \omega^k = 0 \ \forall k, \ 1 \leq k \leq n-1$ (\cite{C}, \cite{X}, \cite{FT})

\item $\omega^l \wedge \ddb \omega^k = 0$ ($(l / k)-$SKT, see \cite{IP}, in particular $\omega^{n-1-k} \wedge \ddb \omega^k = 0$, called {\it generalized $k-$Gauduchon} condition).
\end{enumerate}

In the papers cited in \cite{IP}, examples are given to compare these classes of manifolds; moreover, several applications to physics and geometry are indicated.

As regards our point of view, obviously $k-$SKT implies $k$PL, but we cannot give a characterization by means of positive currents: this motivates our choice.

In the same vein, other kinds of classes of manifolds have been considered: for instance, those where every Gauduchon metric is a strongly Gauduchon metric (in our notation, every $(n-1)$PL metric is also a $(n-1)$S metric), see \cite{PU}, or every $(n-1)$S metric is also an $(n-1)$WK metric, see \cite{FY}. There are important ties with Aeppli cohomology, but this kind of problems is not in the spirit of the present paper.

 \bigskip
 Nevertheless, the link between $\omega$ and $\omega^p$ is important, to answer the following natural question:
 
 Is a \lq\lq$1-$K\"ahler\rq\rq manifold also  \lq\lq$p-$K\"ahler\rq\rq, $\forall  p \geq 1$?
 
 As we said, this is obvious in case K, not in case WK or PL.
 
 In case pS, if we consider the conditions given in Remark 2.4, we are in the same troubles as before. But when the question is translated as follows: \lq\lq Is a 1S manifold also a pS manifold?\rq\rq, the answer is positive.
 
 Indeed, let $\psi$ be a closed 2-form, whose $(1,1)-$component is $\omega > 0$. Then $\psi^p$ is closed too, and its $(p,p)-$component is given by $\omega^p + \zeta$, where $\zeta$ is a sum of $(p,p)-$forms of this kind: $\omega^k \wedge (\sigma_{p-k} \eta \wedge \overline{\eta})$, $\eta \in {\E}^{p-k,0}$, hence 
 $\zeta$ is a positive form, so that $\omega^p + \zeta > 0.$
 \medskip

In general,  we can consider the following classes of manifolds (of dimension $n \geq 3$), for $1 \leq p \leq n-2$:
\medskip

Class \lq\lq $(*)^p$\rq\rq : When $N$ is \lq\lq $p-$K\"ahler\rq\rq, then it is \lq\lq $q-$K\"ahler\rq\rq \ $\forall  q \geq p$.
\medskip

The results we proved give:

\begin{thm} Let $M$ and $N$ be compact complex manifolds, with ${\rm dim} N = n > m = {\rm dim} M$, and let $f : N \to M$ be a holomorphic submersion. Suppose $M$ is K\"ahler. Then:

(i) $N$ belongs to class \lq\lq $(*)^p$\rq\rq,  $\forall  \ p \geq n-m$.

(ii) If $n > m+1$, condition \lq\lq $(HC)$\rq\rq implies that $N$ belongs also to \lq\lq $(*)^{n-m-1}$\rq\rq.
\end{thm}

{\bf Proof.} The first assertion is nothing but Theorem 3.6(1), the second one comes from Theorem 3.6(1) and (3).

\medskip
For instance, $M_{u,v}$ cannot have any degree of K\"ahlerianity, since it is not balanced.
\medskip

As we said before, every compact manifold belongs to $(*)^1_K$ and $(*)^1_S$. On the contrary, conditions $(*)^1_{WK}$ and $(*)^1_{PL}$ are satisfied at least when $N$ has a holomorphic submersion with $1-$dimensional fibres on a K\"ahler basis.

In this situation, the sole condition (\lq\lq HC\rq\rq) implies that $N$ is \lq\lq $p-$K\"ahler\rq\rq, $\forall  p \geq 1$.

\bigskip
\bigskip

\end{document}